\documentclass[12pt]{article}       
\usepackage{amsfonts,amsmath,amssymb}
\newtheorem{Th}{Theorem}
\newtheorem{cor}[Th]{Corollary}
\newtheorem{lem}[Th]{Lemma}

\newenvironment{Pf}
{\noindent\textbf{Proof.}\ \ }{\hfill $\Box$ \medskip}

\title{Online List Colorings with the Fixed Number of Colors}

\begin{document}

\maketitle

\begin{center}{\bf Abstract}\end{center}
\indent\indent
The online list coloring is a widely studied topic in graph theory. 
A graph $G$ is 2-paintable if we  always have a strategy to complete 
a coloring in an online list coloring of $G$ in which each vertex has 
a  color list of size 2.  
In this paper, we focus on the online list coloring game in which 
the number of colors is known in advance. 
We say that $G$ is $[2,t]$-paintable if 
we always have a strategy to complete a coloring in an online list coloring of $G$ 
in which we know that there are exactly $t$ colors in advance, and 
each vertex has a color list of size 2.

Let $M(G)$ denote the maximum $t$ in which $G$ is not $[2,t]$-paintable,   
and $m(G)$ denote the minimum $t \geq 2$ in which $G$ is not $[2,t]$-paintable.   
We show that if $G$ is not 2-paintable, then $2 \leq m(G) \leq 4,$ and 
$n \leq M(G) \leq 2n-3.$ 
Furthermore, we characterize $G$ with $m(G)\in \{2,3,4\}$ and 
$M(G) \in \{n, n+1, 2n-3\},$  respectively.

\section {Introduction}
 
The concept of list coloring was introduced by Vizing \cite{Vizing} 
and by Erd\H os, Rubin, and Taylor \cite{Erdos}.
For each vertex $v$ in a graph $G$, let $L(v)$ denote a list of 
colors available at $v$.
A \emph{  $k$-list assignment} $L$ of a graph $G$ is a list assignment 
$L$ such that $|L(v)|=k$ for each $v\in V(G)$.
A proper coloring $c$ such that $c(v)\in L(v)$ for each $v\in V(G)$ 
is said to be a \emph{list coloring}. 

Schauz \cite{Sch} and Zhu \cite{Zhu} independently introduced an online variation of list coloring. 
On each round $i,$ \emph{Painter} sees only the \emph{marked set} $V_i$ of vertices whose lists contain the color $i.$ 
Painter has to choose an independent subset $X_i$ of $V_i$ to get the color $i.$ 
In the worst case, it can be viewed in the game setting that an adversary, called \emph{Lister}, 
chooses $V_i$ on each round $i$ to prevent a coloring.

Let $f$ be a function from $V(G)$ to the set of nonnegative integers.  
We say that $G$ is $f$\emph{-paintable} if Painter can guarantee coloring all vertices 
with no vertex $v$ is marked more than $f(v)$ times. 
It can be viewed that $f(v)$ is the number of colors that is contained in the list of $v.$  
We write $f \cong k$  if $f(v) = k$ for each vertex $v.$ 
When $G$ is $f$-paintable and $f \cong k$, we say that $G$ is $k$-paintable.

In this paper, we let $(G,f)$ denote the game on a graph $G$ 
with $f$ as the aforementioned function. 
The game $(G,f)$ \emph{contains} $(H,h)$ means that 
$H$ is a subgraph of $G$ and $h(v)=f(v)$  for each $v \in V(H).$ 
 
Three particular functions $f',f^*$ and $f''$ are defined as follows. 
The game $(G,f')$ has $f'(v)=2$ for each $v$ 
except one vertex $u$ which has $f'(u)=1.$ 
The game $(G,f^*)$ is similar to $(G,f')$ except if there is 
a unique vertex $u$ with degree 1, then we always assign $f'(u)=1$ for this vertex.  
The game $(P_n,f'')$ is played on a path $P_n$ for $n \geq 2$ with 
$f''(v)=2$ for each internal vertex $v$ and $f''(u)=1$ for each endpoint $u.$

In this paper, we focus on an online list coloring 
with the given number of rounds to play, 
or equivalently, the given number of colors that appear in all lists. 
Note that the level of information about the number of rounds (colors) 
plays important role for outlining a strategy. 

In this version, Painter knows the number of colors in advance. 
It is reasonable to assume that
Painter knows the number of colors in some applications. 
One maybe more interested in the ``worst case version" of game for Painter, 
namely, Painter does not know the number of colors and Lister knows  
that  Painter does not know the number of colors.  
It is certain that the study of the worst case version is more complicated. 
Nonetheless, the knowledge from the study on this version 
is possibly useful  for facilitating the understanding of other variations.

We say that $G$ is $[f,t]$\emph{-paintable} 
if  Painter guarantees to win in $(G,f)$ with exactly $t$ rounds. 
If $G$ is $[f,t]$-paintable and $f \cong k,$  
then we call $G$ is $[k,t]$\emph{-paintable}.  
Let $M(G,f)$ denote the maximum $t$ in which $G$ is not $[f,t]$-paintable,   
and let $m(G,f)$ denote the minimum $t \geq \max \{f(v): v \in V(G)\}$ 
in which $G$ is not $[f,t]$-paintable.  
If no confusion arises, we may write $M(G)$ and $m(G)$ instead of 
$M(G,f)$ and $m(G,f)$ for $f \cong 2.$ 

The \emph{remaining game} $(G_i,f_i)$ after round $i$ (where  $V_i$ and 
$X_i$ are chosen) is defined recursively as follows. 
Let $(G_0,f_0)=(G,f).$ For $i\geq 1, (G_i,f_i)=(G,f_i)$ where 
$f_i(v)=f_{i-1}(v)-1$ if $v\in V_i,$ and  $f_i(v)=f_{i-1}(v)$ if $v\notin V_i.$ 
If a vertex $v$ is in $X_j$ for some $j\leq i-1,$ 
then we regard $v$ to be colored already and $v$ needs no coloring 
furthermore  in $(G_i,f_i).$

Let $\theta_{p_{1},p_{2},\ldots,p_{r}}$ denote a graph obtained by
identifying all beginnings and identifying all endpoints of $r$ disjoint paths 
having $p_{1},p_{2},\ldots,p_{r}$ edges respectively.
A path $P_{m}$ and a cycle $C_{n}$ intersect at one endpoint of $P_m$ is 
denoted by $P_{m}\cdot C_{n}$. 
Two vertex disjoint cycles $C_{m}$ and $C_{n}$ connected by 
a path $P_{k}$ is denoted by $C_{m}\cdot P_{k}\cdot C_{n}$. 
We always allow $P_k$ in the notation to be $P_1.$   
The \emph{ core} of a graph $G$ is the subgraph of $G$ obtained by 
the iterated removal of all vertices of degree 1 from $G$.

Let $\mathfrak{F}_1 = \{C_{2n+1}\},$  
$\mathfrak{F}_2 = \{C_{m}\cdot P_{k}\cdot C_{n}\}$,  
$\mathfrak{F}_3 = \{\theta_{p,q,r}$ that is not isomorphic to $\theta_{2,2,2n}\}$,  
$\mathfrak{F}_4 = \{\theta_{2,2,2n}$ that is not isomorphic to $\theta_{2,2,2}\}$,  
$\mathfrak{F}_5 =\{K_{2,n}$ where $n\geq 4 \}$, 
and $\mathfrak{F}=\bigcup_{i=1}^5\mathfrak{F}_i.$
\section{Preliminaries and Tools} 

\begin{lem}\label{bip} Assume that $G$ is not 2-paintable. 
A graph $G$ is bipartite if and only if  $m(G)\geq 3.$  \end{lem} 

\begin{Pf} Note that for a game $(G,f\cong 2)$ with exactly two rounds, 
we have $V_1=V_2 = V(G).$  

\emph{Necessity}.  
Assume $G$ is a bipartite graph with partite sets $A$ and $B.$ 
Since $V_1=V_2=V(G),$ Painter can choose $X_1=A$ and $X_2=B$ 
to complete a coloring.  
Thus $m(G)\geq 3.$ 

\emph{Sufficiency}.  
Let $m(G) \geq 3.$ 
In a game of two rounds, 
Painter can choose $X_1$ and $X_2$ which are independent sets 
to complete a coloring.  
Thus $G$ is a bipartite graph with partite sets $X_1$ and $X_2.$ \end{Pf}

\begin{lem}\label{union1} 
Let $G$ be a disjoint union of graphs $H$ and $M.$ 
Let $f(v)=h(v)$ for $v \in V(H),$ and 
$f(v)=m(v)$ for $v \in V(M).$ 
If $H$ is $h$-paintable 
and $M$ is $m$-paintable, then $G$  is $f$-paintable. 
\end{lem}

\begin{Pf}
We prove by induction on the number of uncolored vertices of $G.$ 
Obviously, $G$ is $f$-paintable if $G$ has no uncolored vertices. 
For the induction step, assume that Lister chooses $V_1$ 
in the first round. 
If $V_1 \cap  V(H) \neq \emptyset,$ 
then Painter chooses $X'_1$ that can counter $V_1 \cap  V(H)$ in $(H,h),$ 
otherwise Painter chooses $X'_1 \neq \emptyset.$ 
A set $X''_2 \subseteq V_1 \cap V(M)$ is chosen similarly. 
For $(G,f),$ Painter chooses $X_1 = X'_1 \cup X''_1$ 
to respond for $V_1.$  
The graph $G_1$ in the remaining game $(G_1, f_1)$ 
is the disjoint union of two games that Painter can win. 
Moreover, $(G_1,f_1)$ has fewer uncolored vertices than $(G,f).$  
By induction hypothesis, $(G_1,f_1)$ is $f_1$-paintable. 
Thus $G$ is $f$-paintable. 
\end{Pf}

In a digraph $G$, a set of vertices $U$ is $\emph{kernel}$ of $V'\subseteq V(G)$ 
if $U$ is an independent dominating set of $G[V'].$ 

\begin{lem}\label{tree} 
If $T$ is a tree, then $T$ is $f'$-paintable. \end{lem}
\begin{Pf} 
Let $T$ be an $n$-vertex tree. 
It is clear that Painter wins when $n=1.$ 
Consider $n \geq 2.$ 
Let $u$ be a unique vertex with $f'(u)=1.$ 
Orient $T$ into a digraph in which 
every vertex has in-degree 1 except $u$ which has in-degree 0. 
In the first round, Painter chooses a kernel $X_1$ in $V_1.$ 
Now, $G_1$ in the remaining game $(G_1,f_1)$ is a forest in which 
each nontrivial tree  has all of its vertex $v$ 
satisfying $f_1(v)=2$ except at most one vertex $w$ with $f_1(w)=1.$ 
By induction hypothesis and Lemma \ref{union1}, 
we have $G_2$ is $f_1$-paintable. 
\end{Pf}

\begin{Th} \label{cycleTh}
An odd cycle $C_n$ is not $[2,t]$-paintable if and only if $2 \leq t \leq n.$ 
\end{Th}

\begin{Pf} Consider a game $(C_n,f\cong 2)$ with exactly $t$ rounds. 

\emph{Necessity.} Assume $t \geq n+1.$ Then (i) $V_1 \neq V(C_n),$ 
or (ii) $V_1=V(C_n), |V_i| =1$ for $2 \leq i \leq t=n+1,$ 
and $V_2 \cup \cdots \cup V_{t} =V(C_n).$   
For (ii), Painter just greedily colors a vertex in $V_i$ to win. 

For (i),  $V_1$ induces a union of disjoint paths. 
Orient $V(C_{n})$ to be a directed cycle. 
In the first round, Painter chooses a kernel $X_1$ in $V_1.$ 
Now, the set of uncolored vertices in $(G_1,f_1)$ induces a union of paths in which 
each nontrivial path  has all of its vertex $v$ 
satisfying $f_1(v)=2$ except at most one vertex $u$ with $f_1(u)=1.$ 
By Lemmas \ref{union1} and \ref{tree}, 
Painter has a winning strategy for the remaining game. 

\emph{Sufficiency.} Assume $2 \leq t \leq n.$ 
Lister chooses $V_1 = V(C_{n}).$ 
Regardless of $X_1,$ the remaining game $(G_1,f_1 \cong 1)$ 
has  two adjacent vertices $u$ and $v$ which are uncolored.  
For $t=2,\ldots, t-1,$ Lister chooses $V_i$ to be a set of one vertex other than $u$ and $v.$ 
Finally, in round $t,$ Lister chooses $V_t$ to contain each vertex $w$ with $f_{t-1}(w)=1$ 
(including $u$ and $v$). 
The remaining game $(G_t,f_t \cong 0)$ has $u$ or $v$ uncolored.  
Thus $C_n$ is not $[2,t]$-paintable for $2 \leq t \leq n.$  
\end{Pf}

\begin{lem}\label{subgraph} 
Let the game $(G,f)$ contains $(H,h).$ 
Let $K = \max\{f(v): v \in V(G)-V(H)\}.$ 
If $H$ is not $[h, t]$-paintable, then  $G$  is not $[f,k]$-paintable for 
$\max \{t,K\} \leq k \leq t+ \sum_{v \in V(G)-V(H)}f(v).$ 
In particular, 
 $m(G,f) \leq \max\{K,m(H,h)\}$  
and $M(G,f) \geq M(H,h)+ \sum_{v \in V(G)-V(H)}f(v).$ 
\end{lem}

\begin{Pf}
Lister can win $(G,f)$ with $\max\{t,K\}$ rounds 
by using the strategy similar to one for 
$(H,h)$ with $t$ rounds, except that 
Lister also includes each vertex $v \in V(G)-V(H)$  
in $V_i$ for $i=1,\ldots,f(v).$ 

For $\max\{t,K\}+1 \leq k \leq \sum_{v \in V(G)-V(H)}f(v),$  
Lister has a winning strategy obtained from the above by moving 
vertices in $V(G)-V(H)$ to $V_i$ for $i=\max\{t,K\}+1,\ldots,k$ as needed. 
The remaining follows immediately. 
\end{Pf}

\begin{lem}\label{2-path}  $m(P_2,f'') =1$ and 
$m(P_n,f'') =2$ for $n\geq 3.$ 
\end{lem} 
\begin{Pf} 
The result for  $m(P_2,f'')$ is obvious. 
Consider $n \geq 3.$ 
Let $V(P_n) = \{v_1, \ldots, v_n\}.$ 
If $n$ is even, then Lister chooses $V_1 =  V(P_n) - \{v_1, v_n\}.$ 
The remaining game $(G_1,f_1)$ always has adjacent vertices $u$ and $v$ 
with $f_1(u)=f_1(v)=1.$ 
Lister chooses $V_2 = V(P_n)$ to win the game. 
If $n$ is odd, then Lister chooses $V_1 =  V(P_n) - \{v_n\}.$ 
The remaining game $(G_1,f_1)$ always has adjacent vertices $u$ and $v$ 
with $f_1(u)=f_1(v)=1.$ 
Lister chooses $V_2 =  V(P_n) - \{v_1\}$ to win the game. 
\end{Pf}

\begin{lem} \label{bip1}
If $G$ is a connected bipartite graph with a cycle, then $m(G,f') =3.$ 
\end{lem} 
\begin{Pf} 
Let $u$ be a unique vertex with $f'(u) =1$ 
in a connected bipartite graph $G$ with a cycle $C.$ 

Consider a game $(G,f')$ with two rounds. 
Let $A$ and $B$ be partite sets of $G$ such that $u \in A.$  
Note that  $V_1 = V(G)$ and $V_2=V(G) -\{u\},$or $V_1= V(G) -\{u\}$ 
or $V_2=V(G).$ 
If $u \in V_1,$ then Painter chooses $X_1= A$ and $X_2=B,$  
otherwise Painter chooses $X_1= B$ and $X_2=A.$   
This makes Painter wins. 
Thus $m(G,f') \geq 3.$ 

Next, we show that $m(G,f') \leq 3.$ 
Let $v$ be a vertex in $C$ which is nearest to $u.$ 
Note that $u$ and $v$ can be the same vertex. 
Lister chooses $V_1 =\{x,y\}$ where $xy$ is an edge in $C -\{v\}.$ 
Whatever $X_1$ is, the remaining game $(G_1,f_1)$ 
contains $(P_n,{f''})$ for some $n\geq 2.$ 
The remaining game is not $[f_1,2]$-paintable 
by Lemmas \ref{subgraph} and \ref{2-path}.  
Thus $m(G,f') \leq 3.$ This completes the proof.  
\end{Pf}


\section{Finding $m(G)$}

\begin{lem}\label{cycle2}
If $G$ is bipartite and contains $H \in \mathfrak{F}_2,$ 
then $m(G)=3.$ \end{lem}
\begin{Pf}
Lemma \ref{bip} yields $m(G)\geq 3.$ 
Using Lemma \ref{subgraph}, we only need to show that 
$C_{m}\cdot P_{k}\cdot C_{n}$ 
is not $[2,3]$-paintable  to show $m(G) \leq 3.$ 
First, Lister chooses $V_1 =\{v_1,v_2,w_1,w_2\},$ 
where $v_1v_2$ is an edge in $C_n$,  
$w_1w_2$ is an edge in $C_m,$ and each vertex in $V_1$ is not a cut vertex. 
Regardless of $X_1,$ the remaining game $(G_1,f_1)$ contains 
 $(P_j,{f''})$ for some $j\geq 3.$  
The remaining game is not $[f_1,2]$-paintable 
by Lemmas \ref{subgraph} and \ref{2-path}.  
Thus $m(G)\leq 3.$ 
Lemma \ref{bip} yields $m(G)\geq 3$ which completes the proof. 
\end{Pf}

\begin{lem} \label{theta1} 
If $G$ is bipartite and contains $H \in \mathfrak{F}_3,$  
then $m(G)=3.$ 
\end{lem}
\begin{Pf}
Using Lemma \ref{subgraph}, we only need to show that $\theta_{p,q,r}$ 
where $p,q \geq 3,$ is not 
$[2,3]$-paintable to show $m(G) \leq 3.$ 
Let $P=uw_{1}\ldots w_{p-1}v$, $Q=ux_{1}x_{2}\ldots x_{q-1}v$, 
and $R=uy_{1}y_{2}\ldots y_{r-1}v$ be paths in $\theta_{p,q,r}.$ 
First, choose $V_1 =\{w_1,w_2,x_1,x_2\}.$ 
Regardless of $X_1,$  the remaining game $(G_1,f_1)$ contains 
 $(P_n,{f''})$ for some $n\geq 3.$  
The remaining game is not $[f_1,2]$-paintable 
by Lemmas \ref{subgraph} and \ref{2-path}.  
Thus $m(G)\leq 3.$ 
Lemma \ref{bip} yields $m(G)\geq 3$ which completes the proof. 
\end{Pf}

\begin{lem} \label{theta3}   
If $G$ is bipartite and  contains $H \in \mathfrak{F}_4,$      
then $m(G)=3.$ 
\end{lem}
\begin{Pf}
Using Lemma \ref{subgraph}, we only need to show that $\theta_{2,2,n}$ 
where $n \geq 4,$ is not 
$[2,3]$-paintable to show $m(G) \leq 3.$ 
Let $P=uav$, $Q=ubv$, 
and $R=ux_{1}x_{2}\ldots x_{n-1}v$ be paths in $\theta_{2,2,n}$. 
In the first round, Lister chooses $V_1 = \{u,v,a,x_1,x_{n-1}\}.$ 
Regardless of $X_1,$  the remaining game $(G_1,f_1)$ contains 
the game of $(P_k,{f''})$ for some $k \geq 3.$  
The remaining game is not $[f_1,2]$-paintable 
by Lemmas \ref{subgraph} and \ref{2-path}.  
Thus $m(G)\leq 3.$ 
Lemma \ref{bip} yields $m(G)\geq 3$ which completes the proof. 
\end{Pf}

\begin{lem}\label{theta2} 
If $G \in \mathfrak{F}_5,$  then $m(G)=4.$ 
\end{lem}
\begin{Pf} 
Let partite sets of $G$ be $X= \{x_i : i=1,2,\ldots, n\}$  and $Y= \{u,v\}.$ 
It is well known in the topic of list coloring that 
$K_{2,4}$ is not $L$-colorable 
if $L(u)=\{1,2\},L(v)=\{3,4\},L(x_1)=\{1,3\}, 
L(x_2)=\{1,4\},L(x_3)=\{2,3\},$ and $L(x_4)=\{2,4\}.$ 
Thus $m(K_{2,4}) \leq 4.$ 
Lemma \ref{subgraph} yields $m(G) \leq 4.$

The winning strategy of Painter 
in the game of 3 rounds is as follows: 
Painter colors  both $u$ and $v$ immediately 
after the first $V_i$ that contains $u$ and $v,$   
and greedily colors other legal vertices in other rounds. 
It can be seen that each vertex can be colored. 
Thus Painter wins in the game of 3 rounds. This concludes $m(G)=4.$ 

\end{Pf} 


\begin{lem} \label{core} Assume  $(G,f)$ contains $(H,h)$  
and $H$ is a core of $G.$ 

\noindent (a) If $H$ is $(h,t)$-paintable 
and $2\leq f(v) \leq t$ for each $v \in V(G)-V(H),$  
then $G$ is $(f,t)$-paintable. 

\noindent (b) If $H$ is $h$-paintable and $f(v)\geq 2$  
for each $v \in V(G)-V(H),$ 
then $G$ is $f$-paintable. 

\end{lem}
\begin{Pf} 
(a) We outline Painter's winning strategy for $(G, f)$ as follows. 
Let $F$ be the forest obtained from $G -E(H).$ 
Note that each tree $T$ in $F$ contains at most one vertex $u$ in $H.$ 
Suppose in round $i,$ Lister chooses $V_i.$ 
If $V_i \cap V(H) \neq \emptyset,$   
there is  $X(H)_i$ to counter $V_i \cap V(H)$ 
 in a game $(H,h).$ 
Painter views a game in the part of each tree $T$ in $F$ as a game of $(T, g)$ where 
$g(x) =f(x)$ for each $x \in V(T) -V(H)$ and $g(u)=1$ for a unique vertex in $T \cap H$ 
(if exists.)  
For each tree $T$ and round $i,$ 
Painter considers the marked set $V(T)_i$ in the game $(T,g)$ 
as $(V_i \cap V(T) - \{u\}) \cup (X(H)_i\cap V(T).$ 
Since $g(u) =1,$ Painter  chooses $u$ to be in the set $X(T)_i$ in round $i$ 
if and only if $u \in X(H)_i.$ 
 
Since the coloring of vertices in $V(H)$ which depends on 
Painter's strategy in the game of $(H,h)$ is a winning strategy, 
all vertices in $H$ will be colored. 
By Lemma \ref{tree}, all vertices in each $T$ will be colored. 

(b) is an immidiate consequence of (a). 
\end{Pf}

\begin{lem} \label{m4} Suppose $H$ is the core of a graph $G$ 
and $H$ contains a subgraph in  $\mathfrak{F}_5.$ Then 

\noindent (a) $H \in \mathfrak{F}_5$  
or $H$  contains a subgraph  in $\mathfrak{F}_1 \cup \mathfrak{F}_2 
\cup \mathfrak{F}_3 \cup \mathfrak{F}_4,$    

\noindent (b) $m(G)=4$ if and only if $H \in \mathfrak{F}_5.$  
\end{lem} 

\begin{Pf} 
(a) 
Since $H$ is 2-connected, we can create $H$ from $K_{2,m}$ by successively adding ears 
(an ear is an edge or a path through new vertices connecting two existing vertices)  
or closed ears (a closed ear is a cycle with exactly one existing vertex). 
First, we grow $K_{2,m}$ to be the maximal subgraph $K_{2,n}$ in $H.$ 
For this $K_{2,n},$ let $u,v$ be in the same partite set and $a_1,\ldots,a_n$ be in the other.  
If we cannot add more edges from this point, we have $H = K_{2,n}.$ 
If we can add closed ear, then $G$  contains $C_{s}\cdot P_1 \cdot C_{t}.$ 
If the added ear connects $u$ (or $v$) and $a_i$, then $G$  contains $C_{s}\cdot P_1 \cdot C_{t}.$ 
Consider the case that the added ear has the length $q$ connecting $u$ and  $v.$ 
By maximality of $K_{2,n},$ we have $q \neq 2.$  
Thus if $q$ is odd, then $H$ contains an odd cycle, otherwise $H$ contains 
$\theta_{2,2,2t} \in \mathfrak{F}_4.$ 
Consider the case that the added ear connects $a_1$ and $a_2.$ 
Then the path obtained from the ear plus $a_2v$ is an $a_1v$-path of length at least 3. 
This path together with $a_1v$ and  $a_1ua_3v$ form $\theta_{1,3,q}$ where $q \geq 3.$ 
This completes the first part. 

(b) \emph{Necessity.} 
Suppose the core $H$ of a graph $G$ contains a subgraph  in  $\mathfrak{F}_5.$ 
By (a), $H = K_{2,n},$ 
or $G$  contains $\bigcup_{i=1}^4\mathfrak{F}_i.$  
But the latter case implies $m(G)\leq 3$ by Lemmas \ref{bip}, \ref{cycle2}, \ref{theta1}, 
and \ref{theta3}. 
Hence $H = K_{2,n}$ where $n \geq 4.$ 

\emph{Sufficiency} Suppose $H \in \mathfrak{F}_5.$  
Note that $G$ is bipartite. Thus $G$ is [2,2]-paintable by Lemma \ref{bip}. 
Lemma \ref{theta2} yields $H$ is [2,3]-paintable but not [2,4]-paintable. 
Finally, Lemma \ref{core} yields $G$ is [2,3]-paintable and 
Lemma \ref{subgraph} yields $G$ is not [2,4]-paintable. 
Hence $m(G)=4.$ 
\end{Pf}

\begin{Th}\label{2paint} \cite{Zhu}
A graph $G$ is 2-paintable if and only if the core of $G$ is 
$K_1, C_{2n},$ or $K_{2,3}$.
Equivalently, $G$ is not 2-paintable if and only if the core of $G$ contains 
a subgraph in $\mathfrak{F}.$  
\end{Th}

Now we can classify $m(G)$ for each non-2-paintable graph $G$ as follows.  

\begin{Th}\label{m(G)}
Let $G$ be a non-2-paintable graph. Then $m(G) =2,3,$ or $4.$ 
More specifically, we have 

\noindent (a) $m(G)=2$ if and only if $G$ is not bipartite, 

\noindent (b) $m(G)=3$ if and only if $G$ is bipartite and contains a subgraph in 
$\mathfrak{F}_2 \cup \mathfrak{F}_3 \cup \mathfrak{F}_4,$ 
  
\noindent (c) $m(G) = 4$ if and only if $G$ has a core in $\mathfrak{F}_5.$ 
\end{Th} 

\begin{Pf} 
The statement (a) is exactly Lemma \ref{bip}. 
The statement (c) comes from Lemma \ref{m4}. 
 Let $G$ be a non-2-paintable graph with the core $H.$  
By  Theorem \ref{2paint}, $H$ contains a subgraph in $\mathfrak{F}.$   
By (a) and (c), it remains to consider the case that $G$ is bipartite and $H$ is not in $\mathfrak{F}_5.$ 
By Lemma \ref{m4}, $H$ contains a subgraph in 
$\mathfrak{F}_2 \cup \mathfrak{F}_3 \cup \mathfrak{F}_4.$ 
Since $G$ is bipartite, we have $m(G)\geq 3.$ 
By Lemmas \ref{subgraph}, \ref{cycle2}, \ref{theta1}, and \ref{theta2}, we have 
$m(G)\leq 3.$ 
Thus the remaining case satisfies both $m(G)=3$ and $G$ contains a subgraph in 
$\mathfrak{F}_2 \cup \mathfrak{F}_3 \cup \mathfrak{F}_4.$ This completes the proof. 
\end{Pf}


\section{On $M(G)$} 

Note that $\lg n = \log_2 n.$ 

\begin{lem}\label{2-path1} 
For  $n \geq 3,$ 
$(P_n,f'')$ is not $[f'',t]$-paintable if 
$2\leq t  \leq 2n-2- \lg n .$\\

\end{lem} 
\begin{Pf} 
Let  $n \geq 3$  and $V(P_n)=\{v_1,\ldots,v_n\}.$
We show that $(P_n,f'')$ is not $[f'',t]$-paintable for 
$2\leq t \leq 2n-2-\lg n$ by induction. 
From Lemma \ref{2-path}, we know that $(P_n,f'')$ is not $[f'',2]$-paintable. 
Consequently, the desired statement is true for $n=3.$ 

For $n \geq 4,$ Lister begins with 
$V_1 = \{v_{\lfloor n/2 \rfloor}, v_{\lfloor n/2 \rfloor+1}\}.$ 

Consider the case $v_{\lfloor n/2 \rfloor} \notin X_1.$ 
Then the remaining game 
$(G_1,f_1)$ contains $(P_{\lfloor n/2 \rfloor},f'').$  
By induction and Lemma \ref{subgraph}, 
the remaining game is not $[f_1, t]$-paintable for 
$2\leq t \leq 2\lfloor n/2 \rfloor-2- \lg \lfloor n/2 \rfloor  
+ 2\lceil  n/2\rceil-2=2n-3-\lg (2 \lfloor n/2 \rfloor).$ 
Thus the remaining game is not $[f_1, t]$-paintable for 
$2 \leq t \leq  2n-3-\lg n.$ 
Including the first turn, Lister can win $(P_n,f'')$ with $t$ rounds for 
$3\leq t \leq 2n-2-\lg n.$  

Consider the case $v_{\lfloor n/2 \rfloor+1} \notin X_1.$ 
By induction and Lemma \ref{subgraph}, 
the remaining game is not $[f_1, t]$-paintable for 
$2\leq t \leq 2\lceil n/2 \rceil-2- \lg \lceil n/2 \rceil  
+ 2\lfloor  n/2\rfloor-2 = 2n-3-\lg 2\lceil n/2 \rceil.$  
Note that $\lfloor 2n-3-\lg 2\lceil n/2 \rceil \rfloor 
=\lfloor 2n-3-\lg n \rfloor.$   
Since $t$ is an integer, 
the remaining game is not $[f_1, t]$-paintable for 
$2\leq t \leq 2n-3-\lg n.$ 
Including the first turn, Lister can win $(P_n,f'')$ with $t$ rounds for 
$3\leq t \leq 2n-2-\lg n.$  

\end{Pf}


Let $V(P_m)=\{x_1,\ldots,x_m\}, V(C_n)=\{v_1,\ldots,v_n\}, $ and 
$P_m\cdot C_n$ be obtained from $P_m$ and $C_n$ by identifying $v_n$ with $x_1.$ 
Let $(G,f^*)$ have $f^*(x_m) =1$ and $f^*(v)=2$ for each remaining vertex $v.$ 
Note that $m$ is allowed to be 1.

\begin{lem}\label{1-cycle2} 
If $G=P_m\cdot C_n,$  then 

\noindent (a)  $(G,f^*)$ is not $[f^*,2]$-paintable if and only if $n$ is odd,

\noindent (b) for $t \geq 3,$ $(G,f^*)$ is not $[f^*,t]$-paintable if 
$t \leq 2m+2n-4- \lg (m+ \lfloor n/2 \rfloor).$ 
\end{lem} 

\begin{Pf} 
(a) \emph{Necessity.} If $n$ is even, then $G$ is bipartite.  
Thus $G$ is $[f^*,2]$-paintable by Lemma \ref{bip1}.

\emph{Sufficiency.}  For $n$ is odd, Lister chooses $V_1 =V(G).$
Then the remaining game $(G_1,f_1)$ always contains adjacent uncolored vertices $v$ and $w$ 
in $C_n$ 
such that $f_1(v)=f_1(w)=1.$ Next, Lister chooses $V_2=V(G)-\{x_m\}$ to win the game. 
Thus $(G,f^*)$ is not $[f^*,2]$-paintable.

(b) Lister chooses 
$V_1= \{v_{\lfloor n/2 \rfloor}, v_{\lfloor n/2 \rfloor+1}\}.$ 

If $v_{\lfloor n/2 \rfloor} \notin X_1,$ then the remaining game $(G_1,f_1)$ contains 
$(P_{\lfloor n/2 \rfloor +m},f'')$ which is induced by $\{v_1,v_2,\ldots v_{\lfloor n/2 \rfloor}, 
x_1,x_2,\ldots,x_m\}.$ 
By Lemmas \ref{subgraph} and \ref{2-path1}, $(G_1,f_1)$ is not $[f_1,t]$-paintable if 
$2\leq t \leq 2 (\lfloor n/2 \rfloor +m)-2-\lg (m+\lfloor n/2 \rfloor)
+2(\lceil n/2 \rceil-1) -1 =2m+2n-5-\lg (m+ \lfloor n/2 \rfloor ).$  

If $v_{\lfloor n/2 \rfloor+1} \notin X_1,$ then the remaining game $(G_1,f_1)$ contains 
$(P_{\lceil n/2 \rceil +m-1},f'')$ which is induced by  $\{v_{\lfloor n/2 \rfloor+1},
v_{\lfloor n/2 \rfloor+2},\ldots, 
v_n=x_1,x_2,\ldots,x_m\}.$ 
Thus $(G_1,f_1)$ is not $[f_1,t]$-paintable if 
$2\leq t \leq 
2(\lceil n/2 \rceil + m-1) -2- \lg (m+ \lceil n/2 \rceil-1 )
+2\lfloor n/2 \rfloor -1 
=  2m+2n-5- \lg (m+ \lceil n/2 \rceil-1 ).$
Thus $(G_1,f_1)$ is not $[f_1,t]$-paintable if 
$2\leq t \leq 
2m+2n-5- \lg (m+ \lfloor n/2 \rfloor). $ 

Thus, including the first round, we have $(G,f^*)$ is 
not $[f^*,t]$-paintable if 
$t \leq 2m+2n-4- \lg (m+ \lfloor n/2 \rfloor).$ 
\end{Pf}

Note that the bound in Lemma \ref{1-cycle2} is not sharp if $m$ is large.


\begin{Th} \label{M1}
Let $G$ be a non-2-paintable graph with $n$ vertices. 
Then\\
(a) if $G= C_r \cdot P_k \cdot C_s$ with $r,s \geq 4,$   
then $M(G)\geq n+2,$\\ 
(b) if  $G=\theta_{p,q,r}$ and $p\geq 3, q+r \geq 4,$ then $M(G)\geq n+2,$\\ 
(c) if $G=K_{2,4},$ then $M(G)=n+1=7,$\\   
(d) $M(G) \geq n,$ \\
(e) $M(G)=n$ if and only if $G$ is an odd cycle.
\end{Th}

\begin{Pf} 
(a) Consider $G= C_r \cdot P_k \cdot C_s.$ 
Let $V(C_r) = \{v_1,\ldots, v_r\}$ and $v_r$ be identified with an end vertex of $P_k.$ 
Choose $V_1 = \{v_{\lfloor r/2 \rfloor}, v_{\lfloor r/2 \rfloor +1}\}.$ 
If $v_{\lfloor r/2 \rfloor} \notin X_1,$ then the remaining game $(G_1,f_1)$ 
contains $(C_{s}\cdot P_{k+\lfloor r/2 \rfloor},f^*).$ 
By Lemmas \ref{subgraph} and \ref{1-cycle2}, 
$(G_1,f_1)$ is not $[f_1,t]$-paintable for 
$3\leq t \leq 2(k+\lfloor r/2 \rfloor)+2s-4- 
\lg (k+\lfloor r/2 \rfloor+ \lfloor s/2 \rfloor ) 
+2(\lceil r/2 \rceil -1)-1 
=2k+2r+2s-7- \lg (k+\lfloor r/2 \rfloor+ \lfloor s/2 \rfloor ) .$ 

If $v_{\lfloor r/2 \rfloor+1} \notin X_1,$ then the remaining game $(G_1,f_1)$ 
contains $(C_{s}\cdot P_{k+\lceil r/2 \rceil -1},f^*).$ 
By Lemmas \ref{subgraph} and \ref{1-cycle2}, 
$(G_1,f_1)$ is not $[f_1,t]$-paintable for 
$3\leq t \leq 2(k+\lceil r/2 \rceil -1)+2s-4- 
\lg (k+\lceil r/2 \rceil-1 + \lfloor s/2 \rfloor )
+2\lfloor r/2 \rfloor -1
= 2k+2r+2s-7- \lg (k+\lceil r/2 \rceil-1 + \lfloor s/2 \rfloor ).$ 
Thus $(G_1,f_1)$ is not $[f_1,t]$-paintable for  
$ 3 \leq t \leq 2k+2r+2s-7- \lg (k+\lfloor r/2 \rfloor+ \lfloor s/2 \rfloor ).$ 

Note that $|G|=k+r+s-2.$ 
Thus, including the first round, we have 
$M(G) \geq 2k+2r+2s-6- \lg (k+\lfloor r/2 \rfloor+ \lfloor s/2 \rfloor )
= 2n-2 -\lg (k+\lfloor r/2 \rfloor+ \lfloor s/2 \rfloor )
\geq n+2.$ Note that the last inequality comes from $k\geq 1,r \geq 4,$ 
and $s \geq 3.$

(b) Consider $G=\theta_{p,q,r} \in \mathfrak{F_3} \cup \mathfrak{F_4}.$  
Let $P=uw_{1}\ldots w_{p-1}v$, $Q=ux_{1}x_{2}\ldots x_{q-1}v$, 
and $R=uy_{1}y_{2}\ldots y_{r-1}v$ be paths in $\theta_{p,q,r}.$ 
Choose $V_1 = \{w_{\lfloor p/2 \rfloor}, w_{\lfloor p/2 \rfloor +1}\}.$ 
If $w_{\lfloor p/2 \rfloor} \notin X_1,$ then the remaining game $(G_1,f_1)$ 
contains $(C_{q+r}\cdot P_{\lfloor p/2 \rfloor+1},f^*).$ 
By Lemmas \ref{subgraph} and \ref{1-cycle2}, 
$(G_1,f_1)$ is not $[f_1,t]$-paintable for 
$3\leq t \leq 2\lfloor p/2 \rfloor+2+2(q+r)-4 
- \lg (\lfloor p/2 \rfloor+ 1+ \lfloor (q+r)/2 \rfloor )
+ 2(\lceil p/2 \rceil -1)-1 
=2p+2q+2r-5-  \lg (\lfloor p/2 \rfloor+1+ \lfloor (q+r)/2 \rfloor ).$ 

If $w_{\lfloor p/2 \rfloor+1} \notin X_1,$ then the remaining game $(G_1,f_1)$ 
contains $(C_{q+r}\cdot P_{\lceil p/2 \rceil},f^*).$ 
By Lemmas \ref{subgraph} and \ref{1-cycle2}, 
$(G_1,f_1)$ is not $[f_1,t]$-paintable for 
$3\leq t \leq  2\lceil p/2 \rceil +2(q+r)-4 
- \lg (\lceil r/2 \rceil-1+ \lfloor (q+r)/2 \rfloor )
+ 2\lfloor p/2 \rfloor-1 
=2p+2q+2r-5-  \lg (\lceil p/2 \rceil+ \lfloor (q+r)/2 \rfloor.$ 
Thus $(G_1,f_1)$ is not $[f_1,t]$-paintable for 
$3\leq t \leq 2p+2q+2r-5- \lg (\lfloor p/2 \rfloor+1+ \lfloor (q+r)/2 \rfloor).$ 

Note that $|G|=p+q+r-1.$ 
Thus, including the first round, we have 
$M(G) \geq 2p+2q+2r-4- \lg (\lfloor p/2 \rfloor+1+ \lfloor (q+r)/2 \rfloor )
= 2n-2-  \lg (\lfloor p/2+1 \rfloor+ \lfloor (q+r)/2 \rfloor )
\geq n+2.$ Note that the last inequality comes from $p\geq 3$ and $q+r\geq 4.$

(c)   
Let partite sets of $G$ be $X= \{x_1,x_2,x_3,x_4\}$  and $Y= \{a,b\}.$ 
Observe that Lister has to choose $V_1=\{a,x_1,x_2\}$ 
(or the set of vertices inducing $P_3$) to win the game. 
If $x_1,x_2 \notin X_1,$ then $(G_1,f_1)$ contains 
$(P_3,f'')$ which is induced by $\{x_1,b,x_2\}.$ 
By Lemmas \ref{subgraph} and \ref{2-path1}, 
$(G_1,f_1)$ is not $[f_1,t]$-paintable for 
$3\leq t \leq  2(3)-2- \lg 3  + f_1(a)+f_1(x_3)+f_1(x_4)=9 -\lg 3.$ 
Note that $9 -\lg 3  \geq 7.$

If $a \notin X_1,$ then $(G_1,f_1)$ contains 
$(C_4\cdot P_1,f'')$ which is induced by $\{a,b,x_3,x_4\}.$  
By Lemmas \ref{subgraph} and \ref{1-cycle2}, 
$(G_1,f_1)$ is not $[f_1,t]$-paintable for 
$3\leq t \leq  2(1+4)-4-  \lg 3  +f_1(x_1)+f_1(x_2)=8-\lg 3.$ 
Note that $8-\lg 3  \geq 6.$ 

Including the first turn, we have $M(G) \geq 7 =|G|+1.$

(d)  By Theorem \ref{2paint}, the core of $G$ contains 
a subgraph  $H \in \mathfrak{F}.$  
Lemma \ref{subgraph} yields $M(G) \geq M(H) + 2(n - |H|).$ 
From (a), (b), (c), and Theorem \ref{cycleTh}, 
$M(H) \geq |H|.$ 
Thus $M(G) \geq 2n - |H| \geq 2n -n = n.$ 

(e) \emph{Necessity.} 
In the proof of (d), $M(G) =n$ only if $M(H) =|H|=n.$ 
From From (a), (b), (c), and Theorem \ref{cycleTh}, 
$H$ is an odd cycle $C_n.$ 
If $H \neq G,$ then $G$ contains a smaller odd cycle $C_m$ with $m < n.$ 
Using the proof in (d), $M(G) \geq 2n-m \geq n+2$ which is a contradiction.  
Thus $G$ is an odd cycle. 

The Sufficiency part is an immediate consequence of Theorem \ref{cycleTh}. 
\end{Pf}


Assume  $G$ is a non-$f$-paintable graph. 
Let $q(G,f)$ be the minimum value for $\sum (|V_i| -1)$ that Lister guarantees to have 
where each $V_i$ is a set of marked vertices leading to an \emph{uncolorable} vertex 
(that is an uncolored vertex $v$ with $f_j(v)=0$ for some $j$) with a restriction 
that each vertex $v$ is in at most $f(v)$ sets of $V_i$s.  

For example, consider the game $(P_3,f'')$ where $v_1$ and $v_3$ be endpoints of the path 
and $v_2$ be the remaining vertex. 
Suppose Lister chooses $V_1=\{v_1,v_2\}.$ If Painter does not color $v_1,$ 
then $v_1$ becomes an uncolorable vertex. 
But we cannot conclude that $q(G,h)=|V_1|-1=1$ because Painter may color $v_1.$ 
Painter can choose $V_2=\{v_2,v_3\}$ to guarantee an uncolorable vertex in any cases.  
Thus we can conclude that $q(C_3,h)\leq 2 = (|V_1|-1)+(|V_2|-1).$ 
Lister can continue to choose $V_3=\{v_3\}$ but this does not affect the value 
of $\sum (|V_i| -1)$ and an uncolorable vertex is still uncolorable. 
Generally, if $V_1,\ldots, V_k$ guarantee to force an uncolorable vertex, then 
Lister can choose each remaining $V_i$ to be singleton to retain the value  
of $\sum (|V_i| -1).$ 
It can be seen that this process is unnecessary to continue for finding $q(G,f).$ 

Similarly, if $V_1= \{u\},$ then Painter can color $u.$ 
This does not lead to an uncolorable vertex and the value $|V_1|-1=0$ does not 
affect the value of summation. 
Thus we assume that $V_i$ is not a singleton until an uncolorable vertex occurs. 
If $f(v)=2$ for each $v \in G,$ we just write $q(G)$ instead of $q(G,f).$ 
The next Lemma shows the relation of $q(G,f)$ and $M(G,f).$ 
For convenience, we use $\sum f(v)$ instead of $\sum_{v\in V(G)} f(v).$ 

\begin{lem}\label{qgame} 
$M(G,f) = \sum f(v) - q(G,f).$ 
\end{lem} 
\begin{Pf} 
Since Lister can win in a painting game with $M(G,f)$ rounds, 
Lister can make marked sets $V_1,\ldots,V_{M(G,f)}$ to win a game 
in which each vertex $v$ is in exactly $f(v)$ sets of $V_i$s.   
Note that $\sum (|V_i| -1)= \sum f(v) -M(G,f).$ 
Since  $q(G,f)$ is the minimum value 
of $\sum (|V_i| -1)$ leading to an uncolorable vertex, 
we have $q(G,f) \leq \sum f(v) -M(G,f).$ 
Thus  $M(G,f) \leq \sum f(v) - q(G,f).$ 

Next, by definition of $q(G,f),$ 
Lister can make marked sets $V_1,\ldots,V_k$ to force 
an uncolorable vertex with $q(G,f)=\sum_{i=1}^k (|V_i| -1).$  
After that Lister can choose each $V_i$ for $i = k+1,\ldots, k+\sum f_k(v)$ to 
be a singleton to complete the game $(G,f).$  
Since Painter cannot color an uncolorable vertex, Lister wins by this strategy. 
Consider $\sum f_k(v) = \sum f(v) - \sum_{i=1}^k |V_i| 
= \sum f(v) - \sum_{i=1}^k (|V_i|-1)-k 
=  \sum f(v) - q(G,f) -k.$ 
That is Lister can win $(G,f)$ with $k+\sum f_k(v) =\sum f(v) - q(G,f)$ rounds. 
Thus $M(G,f) \geq \sum f(v) - q(G,f).$ 
This completes the proof. 
\end{Pf} 

Lemma \ref{qgame} implies that finding $q(G,f)$ leads to knowing $M(G,f).$ 
If Painter forces an uncolorable vertex after 
choosing $V_1,\ldots,V_k,$ Painter can minimize 
$\sum (|V_i| -1)$ by choosing $V_j$ to be a singleton for each $j\geq k+1.$ 
But a singleton $V_j$ contributes $0$ in $\sum (|V_i| -1).$ 
Thus to find $q(G,f),$ we may stop counting when an uncolorable vertex occurs.

Next we investigate the condition that $q(G,h)=0, 1,2,$ or $3$ 
where each vertex $v$ has $h(v)=1$ or $2.$  

\begin{lem}\label{q0}
No graph $G$ satisfies $q(G,h)=0.$ 
\end{lem} 
\begin{Pf} To achieve $q(G,h)=0,$ each marked set $V_i$ is a singleton. 
All vertices can be colored which is a contradiction. 
\end{Pf}

\begin{lem}\label{q1}
$q(G,h)=1$ if and only if $(G,h)$ contains $(P_2,f'').$ 
\end{lem} 
\begin{Pf} \emph{Necessity.} 
Let $q(G,h)=1.$ 
Then there is a marked set $V_1=\{a,b\}$ forcing an uncolorable vertex. 
If $a$ and $b$  are not adjacent, then Painter can color both vertices. 
If $h(a)=2,$ then Painter can color $b.$ 
In both situations, an uncolorable vertex does not occur which is a contradiction. 
Thus $a$ and $b$ are adjacent with $h(a)=1.$  
Similarly, $h(b)=1.$ 
 Thus $(G,h)$ contains $(P_2,f'').$ 

\emph{Sufficiency.} 
Assume $(G,h)$ contains $(P_2,f'').$ 
By Lemma \ref{q0}, $q(G,f) \geq 1.$ 
It remains to show that $q(G,f)\leq 1.$ 
Choosing $V_1$ that induces $(P_2,f'').$ 
we have $|V_1| -1=1$ and $V_1$ forces an uncolorable vertex. 
This completes the proof. 
\end{Pf}

We say that a set of vertices $A=\{v_1,v_2,\ldots, v_k\}$ in 
$(G,f)$ induces $(H,h)$ if $A$ induces the graph $H$ and 
$f(v_i)=h(v_i)$ for each $i.$ 

\begin{lem}\label{q2} 
$q(G,h)=2$ if and only if $(G,h)$ does not contain  $(P_2,f''),$ 
but contains $(P_3,f'')$, $(P_4,f''),$ 
or $(C_3 \cdot P_1,f').$  
\end{lem} 

\begin{Pf} \emph{Necessity.} 
Let $q(G,h)=2.$ 
If $(G,h)$ contains $(P_2,f''),$ 
then $q(G,h)=1$ by Lemma \ref{q1} which is a contradiction. 
To have $\sum (|V_i|-1)=q(G,h)=2,$ (i) Lister can choose $V_1$ with size 3 
forcing an uncolorable vertex,     
or (ii) Lister can chooses $V_1$ and $V_2,$ each of which has size 2, 
forcing an uncolorable vertex.

Consider (i). 
Since $(G,h)$ does not contain $(P_2,f''),$  Painter can color 
each $v \in V_1$ satisfying $h(v)=1.$   
An uncolorable vertex does not occur.  
Thus situation (i) is impossible. 

Consider (ii). 
Let $V_1 = \{a,b\}.$ 
If $a$ and $b$ are not adjacent, then Painter can color both $a$ and $b.$ 
Then $V_2$ must induce $(P_2,f'')$ to force an uncolorable vertex which is a contradiction. 
Thus $a$ and $b$ are adjacent. 

For $h(a)=1,$ we assume that Painter chooses $X_1=\{a\},$ 
otherwise Lister can choose $V_2$ to be any 2-set  to have $\sum (|V_i|-1)=2$ and 
an uncolorable vertex. 
Consider the remaining game $(G_1,f_1).$ 
Thus  $q(G_1,f_1)=1.$ 
By Lemma \ref{q2}, $(G_1,f_1)$ contains $(P_2,f'').$ 
Since $(G,h)$ does not contain $(P_2,f''),$ 
this $(P_2,f'')$ contains a vertex $b.$ 
Moreover, there is a vertex $c\neq a$ which has $h(c)=1$ and is adjacent to $b.$ 
Since $(G,h)$ does not contain $(P_2,f''),$ we have  $a$ and $c$ are not adjacent. 
Thus $(G,h)$ contains $(P_3,f'')$ induced by $\{a,b,c\}.$ 

Consider the case $h(a)=h(b)=2.$ 
Since $q(G,h)=2,$ the remaining game $(G_1,f_1)$ always has 
$q(G_1,f_1)=1$ regardless of $X_1.$ 
By Lemma \ref{q1}, $(G_1,f_1)$ contains $(P_2,f'').$  
Thus if $a \notin X_1,$ then there is $c,$ an adjacent vertex of $a,$ 
such that $\{a,c\}$ induces $(P_2,f'').$ 
This also implies $h(c)=1.$ 
Similarly, there exists a vertex $d$ which has $h(d)=1$ and is adjacent to $b,$ 
If $c\neq d,$ then $(G,h)$ 
contains $(P_4,f'')$ induced by $\{a,b,c,d\}.$ 
If $c = d,$ then $(G,h)$ 
contains $(C_3,f')$ induced by $\{a,b,c\}.$

\emph{Sufficiency.} Assume $G$ does not contain $(P_2,f'').$ 
Lemmas \ref{q0} and \ref{q1} imply $q(G,h)\geq 2.$ 
It remains to prove $q(G,h)\leq 2.$ 

Suppose $\{v_1,v_2,v_3\}$ induces $(P_3,f'').$  
Then $V_1=\{v_1,v_2\}$ and $V_2=\{v_2,v_3\}$ force an uncolorable vertex.

Suppose $\{v_1,v_2,v_3,v_4\}$ induces $(P_4,f'').$ 
Then Painter chooses $V_1=\{v_2,v_3\}.$ 
If $v_2 \notin X_1,$ then  $V_2=\{v_1,v_2\}$ forces an uncolorable vertex. 
If $v_3 \notin X_1,$ then  $V_2=\{v_3,v_4\}$ forces an uncolorable vertex.

Suppose  $\{v_1,v_2,v_3\}$ induces $(C_3,f')$  where $f'(v_1)=1.$ 
Then Painter chooses $V_1=\{v_2,v_3\}.$ 
If $v_2 \notin X_1,$ then  $V_2=\{v_1,v_2\}$ forces an uncolorable vertex. 
If $v_3 \notin X_1,$ then  $V_2=\{v_1,v_3\}$ forces an uncolorable vertex. 
 
In each case, We have $\sum (|V_i| -1)=2$ and an uncolorable vertex. 
Thus $q(G,h) \leq 2$ which completes the proof.     
\end{Pf} 

\begin{lem} \label{q3} 
$G$ contains $C_3$ if and only if $q(G)=3.$ 
\end{lem}
\begin{Pf} 
\emph{Necessity.} Let $V(C_3)= \{a,b,c\}.$  
Lister chooses $V_1=\{a,b,c\}. $ 
Since Painter can color at most one vertex, 
we may assume $a$ and $b$ are not colored.  
Choosing $V_2 =\{a,b\}$ forces an uncolorable vertex. 
Thus $q(G) \leq 3.$ 
From Lemma \ref{q0}, \ref{q1}, and \ref{q2}, we have $q(G)\geq 3.$ 
Thus the equality holds. 

\emph{Sufficiency.} Consider the choice of $V_1$ that 
makes $\sum(|V_i|-1)=3$ and leads to an uncolorable vertex. 
Since we want $\sum (|V_i|-1)=2,$ we have $|V_1| \leq 4.$ 
If $|V_1|=4,$ then  remaining $V_i$s are singletons. 
Thus $V_1$ must force an uncolorable vertex. 
But  $f(v)=2$ for each vertex $v,$ an uncolorable vertex does not occur.  
Thus $|V_1| \neq 4.$

Consider $V_1 =\{a,b\}.$ 
Assume that  Painter chooses $X_1 =\{a\}.$ 
By Lemma \ref{q2}, the remaining game $(G_1,f_1)$ contains  $(P_3,f'')$, $(P_4,f''),$ 
or $(C_3 ,f').$  
Since $f(v)=2$ for each vertex $v,$ we have $(G_1,f_1)$ contains  $(C_3 \cdot P_1,f')$ 
and $b \in C_3.$ Thus $G$ contains $C_3.$

Consider $V_1 =\{a,b,c\}.$ 
If $a$ is not adjacent to $b,$ then Painter can choose $X_1 =\{a,b\}.$ 
By Lemma \ref{q1}, the remaining game $(G_1,f_1)$ must contain $(P_2,f'').$ 
This is possible only if $c$ is adjacent to a vertex $v$ with $f(v)=1.$ 
But $f(v)=2$ for each vertex $v.$ This is a contradiction. 
Thus each pair of vertices in $V_1$ are adjacent, that is $G$ contains $C_3.$ 
\end{Pf}

\begin{Th} \label{M2} Let $G$ be a non-2-paintable graphs with $n$ vertices. 
Then the followings hold: \\
(a) $M(G) \leq 2n-3$ for each graph $G,$ \\
(b) $M(G) = 2n-3$ if and only if $G$ contains $C_3.$ 
\end{Th}
\begin{Pf} (a) 
Suppose $M(G) \geq 2n-2.$ 
Thus $q(G) \leq 2$ by Lemma \ref{qgame}. 
But this contradicts to Lemmas \ref{q0}, \ref{q1}, and \ref{q2}.

(b) \emph{Necessity.} $M(G) = 2n-3.$
By Lemma \ref{qgame}, $q(G) =3.$ 
Thus $G$ contains $C_3$ by Lemma \ref{q3}.

\emph{Sufficiency.}  Assume that $G$ contains $C_3.$ 
We have $M(C_3)=3$ by Theorem \ref{cycleTh}. 
Using Lemma \ref{subgraph}, we have $M(G) \geq 2n-3.$ 
Combining with (a), we have the desired equality. 
\end{Pf} 


\section{Further Investigation} 

Using Theorems \ref{M1} and \ref{M2}, we have 
the following corollary. 

\begin{cor} 
If an $n$-vertex graph $G$  is not 2-paintable, 
then $n \leq M(G) \leq 2n-3.$
\end{cor} 

Moreover, we characterizes graphs with 
$M(G)=n$ and graphs with $M(G)=2n-3.$ We turn our attention to 
find the characterizations of $G$ with other values of $M(G).$ 

\begin{lem}\label{last} 
If $n$ is even, then $M(C_{n-1}\cdot P_2)=n+1.$ 
\end{lem} 
\begin{Pf}
Let $G=C_{n-1}\cdot P_2.$ 
By Theorem \ref{M1}, we have $M(G) \geq n+1.$ 

Next we show $M(G) \leq n+1.$ 
Suppose Lister can win in a game of $t$ rounds 
where $t \geq n+2.$ Then (i) $V_1 \not\subseteq V(C_{n-1}),$ 
or (ii) $V_1=V(C_{n-1})$ and $|V_i| =1$ for $2 \leq i \leq t=n+2.$  
Note that in (ii), $V_2 \cup \cdots \cup V_{t} =V(C_n).$   
Thus Painter just greedily colors a vertex in $V_i$ to win. 

For (i),  $V_1$ induces a union of disjoint trees. 
Let $v \in V(G)-V(C_{n-1})$ and $u$ be a neighbor of $v.$    
Orient $V(C_{n-1})$ to be a directed cycle and $u \rightarrow v.$  
In the first round, Painter chooses a kernel $X_1$ in $V_1.$ 
Now, the set of uncolored vertices in $(G_1,f_1)$ induces  
a union of trees in which each tree  has all of its vertex $v$ 
satisfying $f_1(v)=2$ except at most one vertex $u$ with $f_1(u)=1.$ 
By Lemmas \ref{union1} and \ref{tree}, 
Painter has a winning strategy for the remaining game. 
Thus $M(G) \leq n+1$ which completes the proof. 
\end{Pf}

\begin{Th}\label{M3} 
$M(G) = n+1$ if and only if $G$ is $K_{2,4}$ or a 4-vertex graph containing $C_3,$ 
or a core of $G$ is an odd cycle $C_{n-1}.$ 
\end{Th} 

\begin{Pf}
\emph{Necessity.} Let $G$ be a non-2-paintable graph with $M(G)=n+1.$ 
By Theorem \ref{2paint}, $G$ has a subgraph $H \in \mathfrak{F}.$ 
Choose such $H$ with the minimum number of edges.  
By Lemma \ref{subgraph} and Theorem \ref{M1}, 
we have $n+1=M(G) \geq M(H)+ 2(n-|H|) \geq |H| +2(n-|H|)=2n-|H|.$ 
Thus $|H| \geq n-1,$ that is $|H|=n$ or $n-1.$

Consider $|H|=n.$ 
Suppose $H$ is not bipartite. 
If $H$ is not an odd cycle, then $G$ contains $H' \in \mathfrak{F}$ 
such that $e(H')<e(H)$ which contradicts to the choice of $H.$ 
Thus $H$ is an odd cycle. 
Moreover $G=H$ since $G$ cannot have an odd cycle smaller than $H.$ 
But then $M(G)=n $ by Theorem \ref{cycleTh} which is a contradiction.  
Thus $H$ is bipartite. 
This implies $H$ is the graph described in (a), (b), or (c) of 
Theorem \ref{M1}.
But if $H$ is a graph in (a) or (b), then $M(G)\geq M(H) \geq n+2.$ 
Thus $H=K_{2,4}.$ 
If $G \neq K_{2,4},$ then $G$ contains $C_3$ which again 
contradicts to the choice of $H.$ 
Thus $G=K_{2,4}.$

Consider $|H|=n-1.$ 
Suppose $H$ is not an odd cycle. 
Then $H$ is the graph described in (a), (b), or (c) of 
Theorem \ref{M1}.
By Lemma \ref{subgraph} and Theorem \ref{M1}, 
$M(G) \geq M(H)+ 2 \geq (|H|+1)+2=n+2$ which is a contradiction. 
Thus $H$ is an odd cycle. 
Moreover $H$ is an induced subgraph of $G,$ 
otherwise $H$ contains a smaller odd cycle which contradicts 
to the choice of $H.$  
If $H=C_3,$ then $G$ is a 4-vertex graph with $C_3.$ 
Consider the case that $H$ is an odd cycle with length at least 5. 
Let $v \in V(G)-V(H).$ 
If $\deg (v) \geq 2,$ then $\deg(v) =2$ and $G = \theta_{2,2, n-3},$ 
otherwise $G$  has an odd cycle smaller than $H,$ a contradiction. 
By Theorem \ref{M1} (b), $M(G) \geq n+2,$ a contradiction. 
Thus $\deg (v) \leq 1.$ 
This implies $H$ is an odd cycle with $n-1$ vertices 
and it is a core of $G.$

\emph{Sufficiency.} If $G=K_{2,4},$ then $M(G)=7=n+1$ by Theorem \ref{M1}. 
If $G$ is 4-vertex graph containing $C_3,$ then 
$M(G)=5=|n+1$ by Theorem \ref{M2}. 
If a core of $G$ be an odd cycle $C_{n-1},$ then $M(G)=n+1$ by  
Lemma \ref{last}. 
\end{Pf}

\section{Remarks and Open Problems} 

Proceeding to characterize $G$ with $M(G)=n+2$ is  
more involved. First, we need to analyze $M(H)$ where  
$H \in \mathfrak{F}$ more deliberately. 
Moreover, one has to consider the case $|H|=n-2$ and 
other cases carefully. 

Meanwhile, the process to characterize $G$ with $M(G)=3$ can be 
applied to the characterization of $G$ with $M(G)=2n-4.$  
First, begin by characterizing $(G,h)$ with $q(G,h)=3,$ 
and then proceed to characterize $G$ with $q(G)=4.$  
However, the process is clumsy because many more cases arise. 

Thus we propose the first problem. 

\noindent \emph{Problem 1}: Find the efficient method 
to characterize $G$ with $M(G)=n+k$ or $M(G)=2n-k$ for each $k.$ 

Assume that we know a graph $G$ has $m(G)=2$ and $M(G)=2n-3.$ 
Is it true that $G$ is not $[2,t]$-paintable for $2 \leq t \leq 2n-3$? 
The answer is yes. 
By Theorems \ref{m(G)} and \ref{M2}, $G$ contains $C_3.$ 
Using Lemma \ref{subgraph}, we have $G$ is not $[2,t]$-paintable 
for $2 \leq t \leq 2n-3.$ 
This motivates us to ask the second problem. 

\noindent \emph{Problem 1}: 
Suppose that $G$ is not either $[f, t_1]$-paintable or 
$[f, t_2]$-paintable where $t_1 < t_2.$ 
Is it true that $G$ is not $[f,t]$-paintable if $t_1 <t< t_2$? 


\end{document}